\documentclass[twocolumn, letterpaper]{article}

\usepackage{amsmath,amssymb,amsfonts}
\usepackage{algorithmic}
\usepackage{graphicx}
\usepackage{textcomp}

\usepackage{tikz}
\usetikzlibrary{shapes,arrows}
\usepackage{nicefrac}
\usepackage{mathtools, cuted}
\usepackage{lipsum}
\usepackage{float}
\usepackage{afterpage}
\usepackage{placeins}
\usepackage{subcaption}
\usepackage{tabularx}
\usepackage{balance}
\usepackage{textcomp}%
\usepackage{booktabs}
\usepackage{multirow}
\usepackage{accents}

\usepackage[style=ieee]{biblatex}
\usetikzlibrary{patterns,shapes,arrows}
\usetikzlibrary{positioning}
\allowdisplaybreaks

\tikzstyle{decision} = [diamond, draw, fill=blue!20,
    text width=7em, text badly centered, inner sep=0pt]
\tikzstyle{decision1} = [diamond, draw, fill=blue!20,
    text width=7em, text badly centered, inner sep=0pt, node distance = 2.5cm]
\tikzstyle{block} = [rectangle, draw, fill=blue!20,
    text width=13em, text centered, rounded corners, minimum height=3em]
\tikzstyle{block2} = [rectangle, draw, fill=blue!20,
    text width=5em, text centered, rounded corners, minimum height=3em, node distance = 4cm]
\tikzstyle{block3} = [rectangle, draw, fill=red!20,
    text width=5em, text centered, rounded corners, minimum height=3em, node distance = 4cm]
\tikzstyle{line} = [draw, very thick, color=black!50, -latex']
\tikzstyle{cloud} = [draw, ellipse,fill=red!20, node distance=2.5cm, minimum height=2em]

\newtheorem{definition}{Definition}

\newcommand{\Nset}{\mathbb{N}}
\newcommand{\Rset}{\mathbb{R}}

\usepackage[nolist]{acronym}
\begin{acronym}
\acro{OCP}{Optimal Control Problem}
\acro{PMP}{Pontryagin's Minimum Principle}
\acro{MP}{Mathematical Programming}
\acro{DC}{Direct Collocation}
\acro{IRM}{Integrated Residual Methods}
\acro{IRR-DC}{Integrated Regularised Residual Direct Collocation}
\acro{HJB}{Hamilton-Jacobi-Bellman}
\acro{TPBVP}{Two-Point Boundary Value Problem}
\acro{EMPC}{Economic Model Predictive Control}
\end{acronym}

\title{\LARGE \bf Singular Arcs in Optimal Control: \\ Closed-loop Implementations without Workarounds}

\author{Nikilesh Ramesh, Ross Drummond, Pablo Rodolfo Baldivieso Monasterios, Yuanbo Nie
\thanks{This work was supported by an EPSRC Doctoral Training Partnerships EP/W524360/1.}
\thanks{The authors are with the School of Electrical and Electronic Engineering, University of Sheffield, S1~3JD Sheffield, U.K. (Emails: {\tt{\{nramesh2,~ross.drummond,~p.baldivieso,~y.nie\}}} \tt{@sheffield.ac.uk}.)}
}

\begin{document}


\maketitle

\begin{abstract}

Singular arcs emerge in the solutions of Optimal Control Problems (OCPs) when the optimal inputs on some finite time intervals cannot be  directly obtained via the optimality conditions. Solving OCPs with singular arcs often requires tailored treatments, suitable for offline trajectory optimization. This approach can become increasingly impractical for online closed-loop implementations, especially for large-scale engineering problems. Recent development of Integrated Residual Methods (IRM) have indicated their suitability for handling singular arcs; the convergence of error measures in IRM automatically suppresses singular arc-induced fluctuations and leads to non-fluctuating solutions more suitable for practical problems. Through several examples, we demonstrate the advantages of solving OCPs with singular arcs using {IRM} under an economic model predictive control framework. In particular, the following observations are made: \emph{i)} IRM does not require special treatment for singular arcs, \emph{ii)} it solves the OCPs reliably with singular arc fluctuation suppressed, and \emph{iii)} the closed-loop results closely match the analytic optimal solutions.

\end{abstract}

\section{Introduction}
\label{sec:Intro}

Consider the following \ac{OCP} expressed in the general Bolza form
\begin{subequations}
\label{eqn:cont_DOP}
\begin{equation}
    \min_{\substack{x(\cdot),u(\cdot)\ \\ t_0,t_f}} \label{eqn:DOPBolzaObjective}
    \begin{array}{l}
    V_M(x(t_0),x(t_f),u(t_0),u(t_f),t_0,t_f) +  \\ 
    \int_{t_0}^{t_f}\ell(x(t),u(t),t) dt
    \end{array}
\end{equation}
subject to
\begin{align}
     f_d(\dot{x}(t),x(t),u(t),t) = 0,\ & \label{eqn:DEs}\\
    g(\dot{x}(t),x(t),\dot{u}(t),u(t),t) \leq 0,\  &\label{eqn:ineqs}\\
    b_E(x(t_0),x(t_f),u(t_0),u(t_f),t_0,t_f) = 0,\ & \label{eqn:terminal_eq}\\
    b_I(x(t_0),x(t_f),u(t_0),u(t_f),t_0,t_f)\leq 0,\ \label{eqn:terminal_ineq}&
\end{align}
\end{subequations} 

 with $x: \mathbb{R} \to \mathbb{R}^n$ the state trajectory, $u: \mathbb{R} \to \mathbb{R}^m$ the control input, $t_0 \in \mathbb{R}$ and $t_f \in \mathbb{R}$ are, respectively, the initial and final times. $V_M(\cdot,\cdot,\cdot,\cdot,\cdot,\cdot)$ is the Mayer cost functional, $\ell(\cdot,\cdot,\cdot)$ is the Lagrange cost functional, $f_d(\cdot,\cdot,\cdot,\cdot)$ characterises the system dynamics, $g(\cdot,\cdot,\cdot,\cdot,\cdot)$ the path constraints, and $b_E(\cdot,\cdot,\cdot,\cdot,\cdot,\cdot)$ and $b_I(\cdot,\cdot,\cdot,\cdot,\cdot,\cdot)$ are the boundary constraints. When the system dynamics only consists of ordinary differential equations, \eqref{eqn:DEs} becomes~$\dot{x}(t)=f_o(x(t),u(t),t)$. When~\eqref{eqn:cont_DOP} is solved offline, the solution is classed as a trajectory optimisation problem~\cite{betts2010practical}; when implemented in closed-loop, the problem class is understood as being \ac{EMPC} \cite{Grune2014}.

Solving \ac{OCP}s, such as those defined by \eqref{eqn:cont_DOP}, appear in many applied decision-making problems, for example, fast-charging lithium-ion batteries ~\cite{XavierMarceloA.2015Lbcc}. However, as the wealth of literature on this topic can attest, solving \ac{OCP}s can be challenging. Not only can they be computationally intractable, but many interesting and non-trivial structures can emerge in the solutions, complicating problems that appear simple at first glance. One noticeable examples are so-called ``\emph{singular arcs}" which, in a loose sense, concern the inability to obtain  optimal control actions in a finite time interval directly from the optimality conditions. The lack of information about the optimal solution in singular arcs makes finding the control  challenging and, as will be shown in this paper, can lead to spurious solutions when using typical numerical schemes. As well as being interesting from a computational point of view, singular arcs also appear in many practical control problems, including aerospace~\cite{Shen1999} and power systems~\cite{ChangChow1998} applications which motivate their analysis from an applied point-of-view as well. 

To  more explicitly state what singular arcs in \ac{OCP}s are, let us introduce \ac{PMP} ~\cite{kirk2004optimal}. The \ac{PMP} provides necessary conditions for the existence of optimal solutions and, at its core, is the Hamiltonian $\mathcal{H}(\cdot,\,\cdot,\,\cdot,\,\cdot)$ defined by
\begin{equation} \label{eq:PMP:hamiltonian}
    \mathcal{H}({x}, {u}, {p}, t) := \ell({x},{u},t) + \Big{<}{p}, f_o({x},{u})\Big{>},
\end{equation}
with ${p}$ known as the co-states of the problem.  \ac{PMP} states that an optimal trajectory $x^*(\cdot)$ and an optimal input $u^*(\cdot)$ must satisfy the following conditions:

   \begin{subequations} \begin{align}
         {\dot{x}}^*(t) = \frac{\partial \mathcal{H}}{\partial {p}}&({x}^*(t), {u}^*(t), {p}^*(t), t), \label{eq:PMP:xdot}\\
        -{\dot{p}}^*(t) = \frac{\partial \mathcal{H}}{\partial {x}}& ({x}^*(t), {u}^*(t), {p}^*(t), t),\label{eq:PMP:pdot}\\
          \mathcal{H}({x}^*(t), {u}^*(t), {p}^*(t), t) & \leq \mathcal{H}({x}^*(t), {u}(t), {p}^*(t), t) , \label{eq:PMP:globalmin}\\
          \begin{split}
        \label{eq:PMP:transversality}
        [\frac{\partial \Phi}{\partial {x}} ({x}^*(t_f), t_f) - &{p}^*(t_f)]^T \delta {x}_f + \\ [\mathcal{H}({x}^*(t), {u}^*(t), {p}^*(t), t) &+ \frac{\partial\Phi}{\partial {x}} ({x}^*(t_f), t_f)] \delta t_f = 0,
        \end{split}
    \end{align} \end{subequations}
 where $\delta {x}_f$ and $\delta t_f$ are terminal state and time variations. 
 \subsection{Singular arcs}
With the \ac{PMP} and Hamiltonian defined for \ac{OCP}s, the concept of singular arcs in \ac{OCP}s can be introduced.
   \begin{definition}[Singular interval]
     If there exists a time interval $[t_1, t_2]$ of finite duration within which the necessary condition given by \ac{PMP} \eqref{eq:PMP:globalmin} yields no direct information about the optimal control in relation to optimal states and co-states, then interval $[t_1, t_2]$ is known as a singular interval. 
 \end{definition}
 \begin{definition}[Singular arc]
     For singular interval $I_s\subset\Rset$, a singular arc is the control $u(t)$ such that $t\in I_s$.
 \end{definition}

 \subsection{Illustrative example: Aly-Chan problem}
\begin{figure}
    \centering
    \includegraphics[width=\columnwidth]{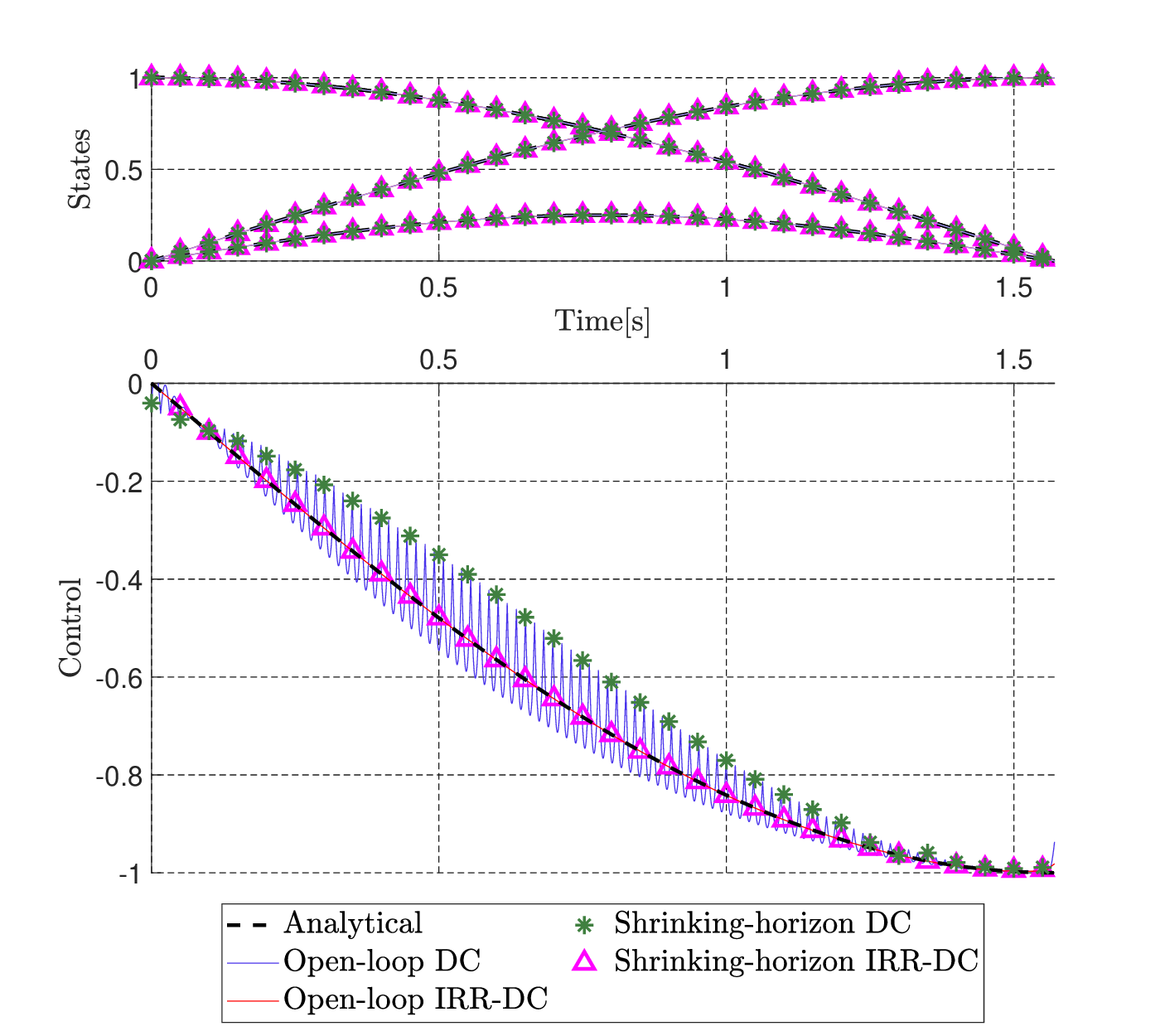}
    \caption{Solutions to the Aly-Chan problem \cite{ALY01041973}. All methods implemented with a fixed mesh of 100 nodes and Hermite-Simpson discretisation. The state plot has $[x_1, x_2,x_3]^\top$. }
    \label{fig:Alychan}
\end{figure}

 To demonstrate the challenges of computing optimal control solutions with singular arcs, consider the classic Aly-Chan problem~\cite{ALY01041973}, which is singular for the entire time interval $[t_0, t_f]$. Figure~\ref{fig:Alychan} compares the performance of different methods that will be introduced later in this paper in both open and closed-loop settings. As the figure shows, the control trajectory from the fixed-mesh open-loop \ac{DC} method (violet solid line) experiences large fluctuations. In contrast, the control trajectory returned by the fixed-mesh open-loop \ac{IRR-DC} method (the red solid line), which is a variant of \ac{IRM}, is much closer to the analytical solution (the black dashed line) derived directly from the PMP optimality conditions. Furthermore, the problem is also implemented in closed-loop as a shrinking-horizon problem with the terminal state being minimized; the closed-loop solutions computed via \ac{DC} (green asterisks) follow a different path to the analytic curve. On the other hand, the closed-loop solution obtained alongside \ac{IRR-DC} matched the analytic curve well.

\subsection{Contributions}

The benefits of \ac{IRM} in reliably solving open-loop \ac{OCP}s with singular arcs have been demonstrated in several studies~\cite{neuenhofen2018dynamic,Nie_2023,nie2025reliablesolutiondynamicoptimization}, outperforming \ac{DC} methods with commonly used fixes to suppress singular arc fluctuations in single-phase setups. This work extends the case-benefit analysis of IRM to closed-loop performance within the popular \ac{EMPC} framework. The primary aim of this paper is to address the following question:  ``\textit{What practical method should be used to numerically solve \ac{OCP}s with singular arcs in closed-loop?}". The following contributions are made in this direction:
\begin{itemize}
    \item Through a review, we highlight that current approaches for addressing singular arcs---both  indirect and  direct methods of \ac{DC}---require prior knowledge of singular arc existence and/or the \ac{OCP} solution structure, making them impractical for real-world engineering problems with high-dimensional \ac{OCP}s. 
    \item We demonstrate the benefits of \ac{IRM} in reliably solving \ac{OCP}s with singular arcs, without additional treatments. The approach is extended to closed-loop implementations and validated against  analytical solutions and analytically derived optimal feedback policies. 
    \item We observe that closed-loop implementations of  DC solutions typically lead to more favourable results: the amplitudes of singular arc fluctuations and the objective value are both reduced in comparison to  open-loop solutions. These results motivate the need to further understand the mechanisms behind these fluctuations.
\end{itemize}

\section{Numerical Methods for Solving Optimal Control Problem}
\label{sec: OptimizationBasedControl}

This section outlines various commonly used numerical methods for solving \ac{OCP}s. As obtaining the analytical solutions of \ac{OCP}s is often intractable in practice, most \ac{OCP}s are solved numerically. Broadly, there are two main solver paradigms for this: the first  are ``\emph{indirect methods}" where the optimality conditions are derived first and then discretized to obtain a numerical solution. The second are ``\emph{direct methods}" where the \ac{OCP} is first discretized via transcription then optimized via \ac{MP}. A detailed survey comparing various variants of direct and indirect methods is available in \cite{rao2009survey}; here, we only provide a high-level overview of the methods. 

\subsection{Indirect method}
The indirect method for solving \ac{OCP}s of the form~\eqref{eqn:cont_DOP} involves finding solutions satisfying certain optimality conditions. A common choice for these conditions, as introduced in Section~\ref{sec:Intro}, are those of the \ac{PMP}. The drawback of using \ac{PMP} is that it only provides necessary conditions for optimality; a trajectory that meets these conditions does not necessarily mean it is the globally optimal trajectory. Nevertheless, \ac{PMP} is still a popular choice for indirect methods because it can reduce the computational overhead. 

Practical implementation of indirect method require analytic expressions for the optimality conditions, and determining these conditions can be hard. Along with the difficulties in initialising the state and co-state variables, the intractability of obtaining these analytic expressions has contributed to the widespread use of the competing direct method for solving \ac{OCP}s in practice.

\subsection{Direct methods} \label{sec:direct:transcription}

In direct methods, the infinite-dimensional continuous-time \ac{OCP}~\eqref{eqn:cont_DOP} is first transcribed into a finite-dimensional \ac{MP} problem, with common choices being quadratic programming and nonlinear programming.

This means that, with direct methods, the trajectories, \emph{i.e.,} the inputs and/or the states, are discretized on a discretization mesh and parameterized by approximating functions $\tilde{x}(\cdot)$ and $\tilde{u}(\cdot)$ with a finite number of decision variables. A finite number of terms then determines the problem, and the remaining challenge of the method is dealing with the integral of~\eqref{eqn:DOPBolzaObjective} and the inequality constraints of~\eqref{eqn:ineqs}. Typically, the integral is evaluated numerically through a quadrature rule while the inequality constraints are only enforced at some selected locations on the discretization mesh~\cite{betts2010practical}.

The main aspect where the different transcription methods differ is how they handle the dynamic constraints of~\eqref{eqn:DEs}. Shooting methods enforce the dynamic relationships via external integrators (\emph{e.g.,}\ differential equation solvers), while direct collocation and integrated residual methods enforce dynamic relationships directly in the MP formulations. In the former case, the feasibility \emph{w.r.t.} the external solver manages the dynamic constraints based on specified error tolerances. In the latter case, posterior error analysis and mesh refinement (MR) are often used to ensure the feasibility of the original \ac{OCP} to specified tolerances.

\section{Solving OCPs with Singular Arcs}
    Although the methods outlined in Section \ref{sec: OptimizationBasedControl} have been found to work for most practical \ac{OCP}s, they can struggle when the solution has singular arcs~\cite{betts2010practical}. Over the years, many classical methods have been adapted to address these issues. The following section gives a brief overview of these adaptions and highlights their suitability for closed-loop control.

\subsection{Overview of commonly used methods}
\subsubsection{Indirect methods for solving singular arcs}\label{sec:indirect}
    Since indirect methods act directly on the conditions of the \ac{PMP}, this optimality condition alone cannot determine the optimal control along singular arcs. A common workaround consists of adding extra conditions to the \ac{OCP} to supply additional information and make the solutions of this augmented problem feasible to compute. This approach is possible for problems where the control: \emph{i)} appears linearly in the Hamiltonian, \emph{ii)} is box-constrained \emph{i.e.,} $u_L \leq u(t) \leq u_U$, and \emph{iii)} does not include any state or mixed state-input constraints \cite{bryson1975applied}. 

    For problems with the structure stated above, switching functions of the form\[S(t) \coloneqq \frac{\partial \mathcal{H}(x(t), p(t),u(t),t)}{\partial u}\]
     define the solution structure via
          \begin{equation} \label{eq:indirect:control:structure}
        u^*(t) = 
        \begin{cases}
            u_L & \text{if } S(t) > 0, \\
            u_U & \text{if } S(t) < 0, \\
            u_S(t) & \text{if } S(t) = 0,
        \end{cases}
    \end{equation}
   with $u_U$ and $u_L$ being upper and lower bounds, and $u_S$ the singular arc.  To be specific, it is noted that at  the switching times $t_s\geq 0$ when $S(t_s) = 0$, the control is defined to change from $u(t_s-\varepsilon) = u_L$ to $u(t_s+\varepsilon) = u_U$, or vice-versa, or a switch from non-singular to singular interval or vice-versa, for some small $\varepsilon >0$. 

    Using the fact that the switching function and all its time derivatives equal zero within the singular interval,  the singular arc can be characterized by analysing higher derivatives of  $S(t)$ until the control $u(t)$ appears explicitly \cite{kirk2004optimal}, \emph{i.e.,}\
    \begin{equation} \label{eq:Singular:SwitchDerivative}
        u(t) = \text{arg}\left\{ \frac{d^{2k}}{dt^{2k}}S(t) = 0\right\}, k\in\Nset.
    \end{equation}
    It is noted that an even number of derivatives is needed to ensure optimality of the singular arc \cite{Powers1980ER} with this condition. Furthermore, for optimal control along singular arcs, the following condition holds:
    \begin{equation} \label{eq:Singular:optimality}
    (-1)^{k} \frac{\partial}{\partial u} \Big{(}\frac{d^{2k}}{dt^{2k}}S \Big{)} \leq 0.
    \end{equation}
   This condition is known as the Kelley optimality condition (or the Generalized Legendre-Clebch condition \cite{KELLEYHENRYJ1964Asvt}),  and has several interesting implications for \ac{OCP}s, as discussed in Section \ref{sec:SMIB}. 

    Solving \ac{OCP}s with singular arcs using the indirect approach leads to several challenges in practice. Firstly, the solution structure and the switching times, \emph{i.e.,} the time $t_s$ at which $ \lim_{t \rightarrow t_s}S(t) = 0$, need to be determined. One way of doing this is to search on a $n-1$-dimensional hyper-surface $\Lambda_s$ such that both singular and non-singular arcs must simultaneously satisfy a junction condition,$$\frac{\partial \mathcal{H}}{\partial u} = 0,$$ with all its time derivatives also equal to zero \cite{Anderson1972}. Applications where this approach has been used can be found in \emph{e.g.,} \cite{ALY01111978}. 

    Secondly, the obtained optimal inputs may not be suitable as feedback policies for closed-loop implementations. Section~\ref{sec:SMIB} illustrates an example where an open-loop stable system becomes unstable when coupled with the derived optimal feedback policy along the singular arc. 
        
    In addition, the complexity of the analytical derivation process can quickly become unmanageable and intractable for real-world problems, even with modern computer-based symbolic operations. As a result of these limitations and challenges, most practical OCPs with singular arcs are solved either using direct methods, or a combination of direct and indirect methods.

\subsubsection{Direct methods for solving singular arcs}

\paragraph{Direct Collocation}
The \ac{DC} method is perhaps the most popular choice for solving a wide range of \ac{OCP}s. With \ac{DC}, the dynamic constraints~\eqref{eqn:DEs} are enforced in the optimization problem as equality constraints based on the weighted residual method of collocation. Effectively, the constraints force the residual error to zero only at some selected locations of the discretisation mesh, known as collocation points. For a time mesh with $Z$ intervals, each with $N^{(z)}$ collocation points $\tau^{(z)}$, the equality constraint implemented in the MP for direct collocation is
\begin{multline*}
f_d\Big(\dot{\tilde{x}}(\tau_i^{(z)}),\tilde{x}(\tau_i^{(z)}),\tilde{u}(\tau_i^{(z)}),\tau_i^{(z)}\Big)=0, \\ \text{for all } i\in\mathbb{I}_{N^{(z)}}, \,z\in\mathbb{I}_Z.
\end{multline*}

When DC is used to solve \ac{OCP}s with singular arcs, the obtained numerical solution may contain large fluctuations, also known as ringing or spurious solutions. This behaviour can be attributed to a combination of factors. Firstly, the switching time for the transition into singular arcs may not be captured exactly by a discretization mesh. The solutions on the singular arc may then need to be adjusted to compensate for the mismatch in switching time. This effect is exaggerated because the sensitivity of input variations to the objective diminishes on the singular arc. Hence, solutions with negligible differences in the objective can have drastically different characteristics, ranging from solutions with no (or small) fluctuations that closely resemble the original continuous-time \ac{OCP} solutions to those with large fluctuations. 

Because of larger residual errors associated with trajectories of high-frequency fluctuations and accounting for the operational limits of actuators, solutions with large fluctuations are generally considered undesirable from an implementation point of view. For example, in the case of fast-charging lithium-ion batteries, these fluctuations could lead to increased cell heat generation and, therefore, to accelerated degradation. Several workarounds have been proposed to suppress or avoid the fluctuations, which are presented below.

\noindent \underline{Control regularization:}
The most widely used strategy to suppress  fluctuating control trajectories along singular arcs is to regularize the control. Typically, this is achieved by including the 2-norm of control actions and/or its rate of change in the cost function using some weighting parameters. Other variations of control regularization has been investigated in~\cite{heidrich2021investigation}.
\noindent \underline{Multiphase approach:}
The time interval is subdivided into segments or phases where different objectives and constraints apply. For \ac{OCP}s with singular arcs, this approach could also be leveraged to enforce the additional analytically derived conditions via the indirect approach on the singular arcs so that the multiphase direct collocation problem is well-determined and the fluctuations are avoided.

\noindent \underline{Hybrid approach:}
A combination of both multiphase and control regularization is used in \emph{Bang-Bang Singular Optimal Control} (BBSOC) \cite{Pager2022}. This approach employs jump functions to detect discontinuities in the control. Singular interval conditions from the \ac{PMP} are then used to classify singular and non-singular regimes. Compared to conventional indirect approaches, the conditional expressions are numerically evaluated. Applying this method involves adding a regularization term $\delta$ to the cost functional of the original \ac{OCP} within the singular regime, given by
    \begin{equation*}
        \delta = \frac{\epsilon}{2} \int_{t_{s_1}}^{t_{s_2}} (u(t) - \gamma_e(t))^2 dt,
    \end{equation*}
    where $e= 1, \,2,\,..., $ is the iteration number, $\epsilon$ is the regularization weighting parameter, and $\gamma_e$ is set to zero for the first iteration and then it  is iteratively updated with cubic polynomial approximation of the control at iteration $e-1$ until $\delta \approx 0$. By regularising only within the singular interval, BBSOC reduces the effect of regularization on the sub-optimality of the resulting solution to the OCP.

\paragraph{Integral Residual Method} \label{sec:IRM}
The main idea behind \ac{IRM} is the following: instead of enforcing the residual to be zero only at selected locations  (as the collocation points in \ac{DC}), \ac{IRM} minimizes and/or constrains the integral of the residual along the \emph{entire} solution trajectory. In practice, this integral can be computed through a $Q^{(z)}$-point quadrature rule inside each mesh interval, as
\begin{equation*}
\sum_{z=1}^{Z}\sum_{i=1}^{Q^{(z)}}w_{i}^{(z)}\begin{Vmatrix}\alpha \circ f_d\Big(\dot{\tilde{x}}(q_i^{(z)}),\tilde{x}(q_i^{(z)}),\tilde{u}(q_i^{(z)}),q_i^{(z)}\Big)\end{Vmatrix}^2_{2},
\end{equation*}
with \emph{quadrature weights} $w_i^{(z)}$, \emph{quadrature abscissae} $q_i^{(z)}$ and $\alpha \in \mathbb{R}^{n_f}$ additional weighting parameters, e.g.\ to account differences in the numerical range of dynamic constraints. 

In comparison to \ac{DC}, \ac{IRM} can: \emph{i)} obtain solutions of higher accuracy for a given mesh, \emph{ii)} facilitate the flexible trade-off between solution accuracy (in terms of dynamics feasibility) and optimality, and \emph{iii)} provide reliable solutions for challenging problems, including those with singular arcs~\cite{Nie_2023}.

The success of the \ac{IRM} method has spawned many versions, including the quadrature penalty method~\cite{neuenhofen2018dynamic, neuenhofen2023numerical}, direct alternating integrated residuals~\cite{Nie_2023, nie2020efficient}, \ac{IRR-DC}~\cite{nie2025reliablesolutiondynamicoptimization}, as well as a multiple shooting variant with a tailored Runge-Kutta integrator for the differential equations~\cite{harzer2025integrationerrorregularizationdirect}. An overview of different IRM variants is available in~\cite{nie2025reliablesolutiondynamicoptimization}. In this paper, we will use the \ac{IRR-DC} approach because it retains the merits of \ac{DC} while improving solution accuracy. \par

\subsection{Suitability analysis for closed-loop implementation}

In contrast to offline trajectory optimization, where multiple design iterations can be performed to mitigate the impact of singular arcs, the closed-loop implementation of \ac{OCP}s presents additional challenges. Specifically, methods requiring special treatment for singular arcs may become impractical in this setting. 

Firstly, it is often not feasible to know the optimal solution structure in advance. Changes in initial conditions during EMPC updates may then lead to changes in the solution structure. This makes methods reliant on analytical derivations, such as the multi-phase approach with DC, unsuitable for closed-loop applications. These methods typically require complicated procedures to identify the solution structure and/or derive conditions related to singular arcs based on that structure. 

Secondly, in practical engineering problems with economic objectives (e.g., minimum-time or minimum-energy), it may not even be possible to know whether singular arcs are present beforehand. While this analysis may seem trivial for simple examples (such as ones shown in Section~\ref{sec: Example}), it becomes intractable as problem complexity increases. 

In control regularisation, the weighting parameters ideally should be chosen so that the penalty on the control input is large enough to specify a preference amongst the solutions. This way, solutions with large control input and input fluctuations are disfavoured. At the same time, the penalty on the control input should not contribute significantly to the cost, leading to a loss of optimality. The success of control regularization depends on the careful management of relative weighting between the original objective and the regularization term. This becomes increasingly difficult in closed-loop systems with many input variables:
\begin{itemize}
    \item The numerical range of the economic objective can vary greatly as the system evolves, making fixed weighting parameters unsuitable. 
    \item Without knowledge of the singular arc, it may be difficult to determine which input variables should be regularized in advance.
\end{itemize}

In practice, a compromise may be necessary to suppress fluctuations along singular arcs reliably in a closed-loop setting, where highly suboptimal solutions might need to be accepted. Iteratively adjusting the regularization, as done in the hybrid approach of BBSOC, requires solving the \ac{OCP} repeatedly, which impacts the computational performance of the \ac{EMPC}.

In contrast to \ac{DC} approaches, which require explicit procedures for \ac{OCP}s with singular arcs, \ac{IRM} directly suppresses solution fluctuations along singular arcs through its error-handling schemes. This provides both flexibility and generality. In \ac{IRM}, different solution candidates on the singular arc with negligible differences in cost can be distinguished by their integrated residual errors: larger fluctuations in the solution typically lead to larger errors along the trajectory when approximated with parameterized functions. As a result, the solution with the smallest fluctuations is often the most accurate according to the \ac{IRM} error metrics. This ability to suppress singular arc fluctuations without requiring additional treatments for each specific \ac{OCP} is why \ac{IRM} is highly suited for online solutions in closed-loop settings. It eliminates the need for posterior analysis and \ac{OCP} design iterations aimed at singular arc mitigation, which would otherwise be inefficient and impractical.

\section{Numerical Examples}  \label{sec: Example}

In this section, we investigate the closed-loop performance of DC and IRR-DC for OCPs with singular arcs, under the framework of EMPC with receding horizon control. To allow for effective comparisons, we select two commonly used examples with simple problem formulations to facilitate the derivation of optimal feedback policies via the indirect approach for reference. The primary goal of 
these examples is to assess the differences in closed-loop performance between DC and IRR-DC for real-world engineering problems, with no prior knowledge of the existence of singular arcs. 

The first example is a classic second-order singular control problem, extensively studied in the literature (e.g., \cite{JacobsonD.1970Coos, pagerThesis, ANDRESMARTINEZ2020106923, ALY01111978}) using both indirect and direct methods. This example serves as a basis for comparison with the solutions derived from IRR-DC. Additionally, it highlights the effects of closed-loop implementation on direct method solutions along the singular arc. In the second example, we explore a practical problem from the power systems domain. Here, we derive the optimal feedback policies and demonstrate that their implementation destabilizes the stable open-loop system. However, when using EMPC with OCPs solved via the direct method of IRR-DC, closed-loop results closely match the analytically derived optimal feedback policy.  All problems are transcribed using the optimal control software \texttt{ICLOCS2}~\cite{ICLOCS2}, and numerically solved to a tolerance of $10^{-9}$ with the NLP solver \texttt{IPOPT}~\cite{wachter2006implementation}.

\subsection{Second Order Singular Control Problem}
 The aim of this example is to compare the performance of IRR-DC to indirect and DC-based direct approaches (such as BBSOC) and to quantify the optimality gap in relation to the analytically found solution. The second order singular control problem \cite{JacobsonD.1970Coos} involves minimising the following objective
\begin{subequations}\begin{equation*}
    J  = \int_{0}^{5}\frac{1}{2}({x_1(t)}^2 + x_2(t)^2) \, dt,
\end{equation*}
subject to 
\begin{equation*}
            \dot{x}_1(t) = x_2(t), \quad
            \dot{x}_2 (t) = u(t), \quad |u(t)| \leq 1. 
    \end{equation*}\end{subequations}

 The optimal control policy on the singular arc can be found analytically as $u^*(t) = x_1(t)$. But, as was highlighted earlier, finding the switching time $t_s$ can be challenging. For example, Aly \cite{ALY01111978} searched for a point on a hyper-surface where both the singular and non-singular arcs satisfied the junction constraints. This additional effort in determining the solution structure is not ideal when dealing with online optimal control of large-scale problems.

.

For methods that do not require prior knowledge of the solution structure, Table \ref{tab:ModAly:compare} compares the  costs and  optimality gaps for the various methods against the analytical optimal solution. In the table, a distinction is made between the objectives reported by the solvers and those practically achieved through numerical simulations.  

\begin{figure}[h]
        \centering
        \begin{subfigure}{\columnwidth}
            \includegraphics[width=\columnwidth]{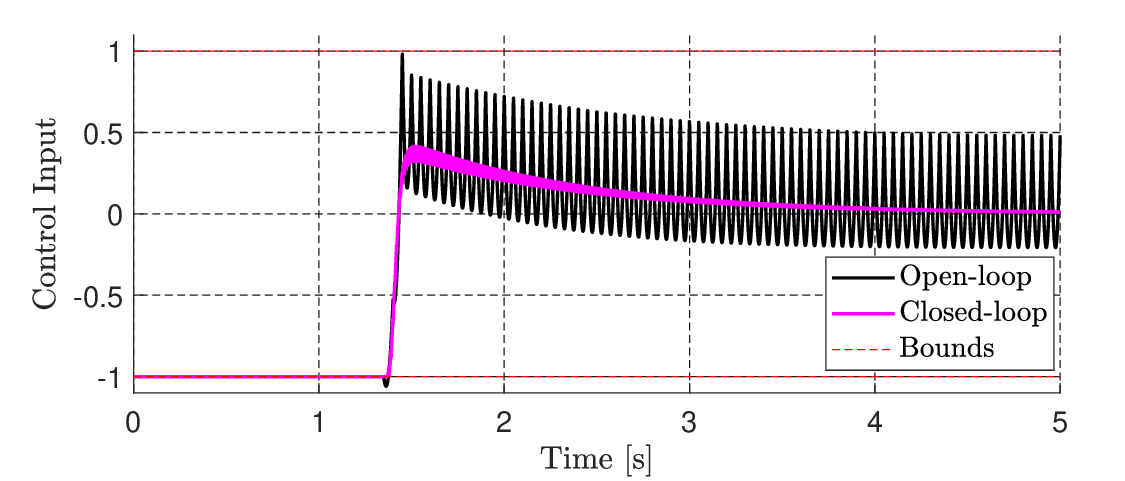}
            \caption{DC}
            \label{fig:ModAly:DC}
        \end{subfigure} \\
        \begin{subfigure}{\columnwidth}
            \includegraphics[width=\columnwidth]{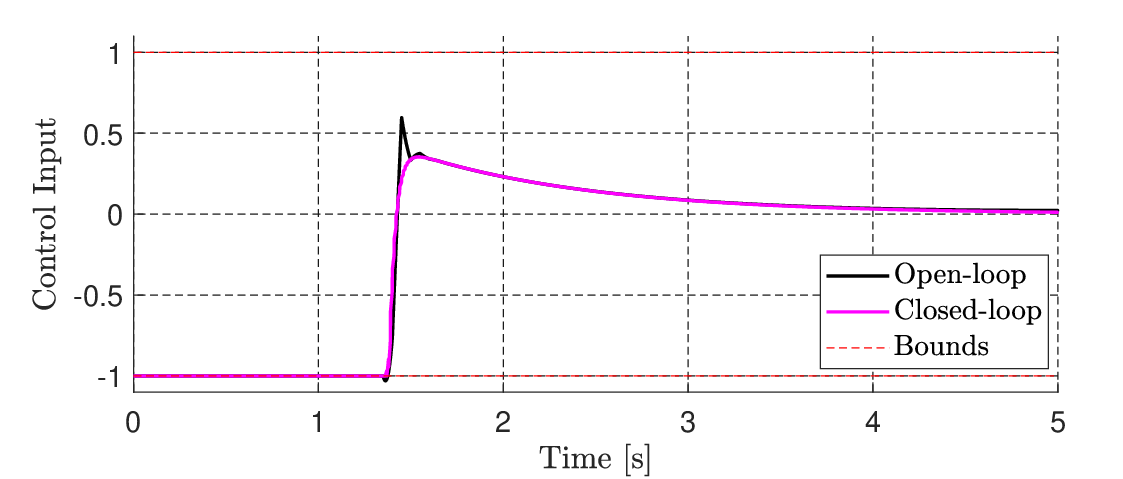}
            \caption{IRR-DC}
            \label{fig:ModAly:IRR}
        \end{subfigure}
        \caption{Comparison of different solutions to the second order singular control problem. Closed-loop implementations use a time step of $0.01s$ and a horizon of $5s$.}
        \label{fig:ModAly}
    \end{figure}

The differences in objective values are negligible when comparing the open-loop \ac{OCP} solutions. However, Figure~\ref{fig:ModAly:DC} shows that the DC solution contains large fluctuations. Based on~\cite{pagerThesis}, BBSOC, using a multi-phase \ac{OCP}, achieves a similar cost while suppressing fluctuations. \ac{IRR-DC} successfully suppressed fluctuations under a single-phase setting without prior knowledge nor special treatment for the singular arc.

\begin{table}[t] 
\begin{minipage}{\linewidth} 
  \small
	\begin{center}
	\caption{Comparison of solutions to the singular control problem via different methods. Optimality gap is computed as $\frac{|J_a - J_m|}{J_a}$, $J_m$ the cost from other methods and $J_a$ the analytical optimal cost reported in the literature.}
	\label{tab:ModAly:compare}
		\begin{tabular}{c|c|c|c}
		 \multicolumn{2}{c|}{\multirow{2}{*}{\textbf{Method}}}   & \textbf{Minimized} & \textbf{Optimality} \\
		 \multicolumn{2}{c|}{} & \textbf{cost} & \textbf{gap}\\
		\hline \multicolumn{2}{c|}{\multirow{2}{*}{Analytic}} & \multirow{2}{*}{$0.37699^{a}$} & \multirow{2}{*}{-}\\
		\multicolumn{2}{c|}{} & &\\
		\hline & \multirow{2}{*}{BBSOC} & \multirow{2}{*}{$0.37699193$\footnote{As reported in \cite{Pager2022}.}} & \multirow{2}{*}{$0.0005\% $} \\
		  & & &\\
          \cline{2-4} OCP &{Open-loop} & \multirow{2}{*}{$0.37699186 $} & \multirow{2}{*}{$0.0005\%$}\\
         Solution & DC& \\
		 \cline{2-4} &{Open-loop} & \multirow{2}{*}{$0.37699222$} & \multirow{2}{*}{$0.0006\% $}\\
       & IRR-DC& \\
       \hline &{Open-loop} & \multirow{2}{*}{$0.37775$} & \multirow{2}{*}{$0.2020\%$}\\
         & DC& \\
		 \cline{2-4} &{Open-loop} & \multirow{2}{*}{$0.37770$} & \multirow{2}{*}{$0.1884\%$}\\
       Simulated & IRR-DC& \\
        \cline{2-4} Solution & {Closed-loop}  & \multirow{2}{*}{$0.37708$} & \multirow{2}{*}{$0.0242\%$}\\
        & DC & \\
        \cline{2-4} & {Closed-loop}  & \multirow{2}{*}{$0.37706$} & \multirow{2}{*}{$0.0198\% $}\\
        &IRR-DC & \\
		\hline
		\end{tabular} 
	\end{center}
    \end{minipage}
\end{table}

In the \ac{OCP} solution, the open-loop \ac{IRR-DC} cost is slightly higher than the \ac{DC} cost. However, when comparing the simulated solutions, implementing the \ac{IRR-DC} solution led to a smaller optimality gap than \ac{DC}. This is because the low objective returned by the \ac{OCP} solver can be misleading, being obtained at the expense of dynamic feasibility. \ac{IRR-DC}'s advantage over DC in implementation is expected to be even larger when practical actuation mechanism rate limits are considered. 

Moving to a closed-loop implementation of the \ac{OCP} solution under a receding horizon setup, the objective for both \ac{DC} and \ac{IRR-DC} saw substantial decreases, showing the benefit of feedback. The \ac{IRR-DC} solution reduced the optimality gap by 18\%  compared to  DC methods when in closed-loop. This reduction in the optimality gap was combined with the input trajectory from the \ac{IRR-DC} method being smooth, in contrast to the fluctuating \ac{DC} one, as evidenced in Figure~\ref{fig:ModAly}. The results demonstrate the competitiveness of \ac{IRR-DC} in closed-loop implementation of \ac{OCP}s with singular arcs. 

Figure~\ref{fig:ModAly:DC} demonstrates the closed-loop implementation of the \ac{DC} solution reduces the amplitude of singular arc fluctuations. This reduction is related to the choice of \ac{EMPC} time steps: a faster feedback loop appears beneficial in alleviating singular arc fluctuations in closed-loop. A detailed investigation of the mechanism behind this behaviour is planned for future work.

\subsection{Single machine on infinite bus} \label{sec:SMIB}
Next, we consider a Single Machine Infinite Bus (SMIB) and the problem of determining the optimal recovery profile in response to a disturbance \cite{ChangChow1998}. The SMIB system has two states ${x} = [x_1, x_2]^\top$, with $x_1$ being the rotor angle and $x_2$ being the frequency relative to the infinite bus. The dynamic model of the system satisfies\cite{jones2024modelpredictivebangbangcontroller}:
\begin{equation} \label{eq:SMIB:xdot}
        \begin{split}
            \dot{x}_1(t) &= x_2(t),\\
            \dot{x}_2(t) &= \frac{P_M - D x_2(t)}{2H} - \frac{P_E}{2H}\sin(x_1(t) + \delta_{ep})u(t),
        \end{split}
    \end{equation}
where $P_M$ and $P_E$ are, respectively, the mechanical  and electrical powers, $D$ is the damping ratio, and $H$ is the inertia constant. The initial conditions of the system are ${x}(t_0) = [1.5,15]^\top$. The input $u(t)$ corresponds to the scaling factor of maximum electrical power controlled via switchable capacitors and is bounded by $|u(t)| < 1$.

The problem is to drive the system back to equilibrium in the wake of disturbances, such as load changes, whilst satisfying input constraints. The cost function is
\[J = \int_{0}^{4} ||{x}(t)||^2dt.\]

A singular analysis is applied to the \ac{OCP} described above to show that the SMIB optimal control has a singular arc.
The solution structure of the \ac{OCP} is then derived using the singular interval conditions and the theory introduced in Section~\ref{sec: OptimizationBasedControl}. Then, following the procedures in Section~\ref{sec:indirect}, the optimal control policy is
\begin{equation*} \label{eq:SMIB:control:structure}
    u(t) = 
        \begin{cases}
            -1 & \text{if } p^*_2(t) < 0, \\
            +1 & \text{if } p^*_2(t) > 0, \\
            u_S & \text{if } p^*_2(t) = 0,
        \end{cases}
    \end{equation*}
with
\begin{equation}
\label{eq:SMIB:control:feedbackpolicy}
    u_S\left(t\right)=-\frac{H\,{C_2 }^2 \,{\left(\frac{2\,x_1 \left(t\right)}{{C_1 }^2 }-\frac{{P_M}-{D}\,x_2 \left(t\right)}{H\,{C_2 }^2 }\right)}}{{P_E}\,\sin \left({\delta_{ep}}+x_1 \left(t\right)\right)},
\end{equation}
where $C_1 = 3,\, C_2 = 30, \,P_M~=~1, \,P_E~=~4.3214,$ $\, D~=~0.03,\, \delta_{ep} = 0.235$. This control law needs to satisfy Kelley's condition \eqref{eq:Singular:optimality} to be optimal. For the SMIB, this condition reduces to
\begin{equation*}
   \sin \left({\delta_{ep}}+x_1 \left(t\right)\right) \geq 0.
\end{equation*} 

\subsubsection{System Stability with the optimal feedback policy}\label{sec:dstable}

Within the singular intervals, $p^*_2(t) = 0$ and $u(t) = u_S(t)$. Substituting these conditions into the dynamic system \eqref{eq:SMIB:xdot} gives a (locally) autonomous system of the form 
\begin{equation} \label{eq:SMIB:autonomous:sys}
   \dot{{x}}(t) = A \, {x}(t) = \begin{bmatrix}
       0 & 1 \\
       \left(\dfrac{C_2}{C_1}\right)^2 & 0
   \end{bmatrix}
   \begin{bmatrix} 
    {x}_1(t) \\
    {x}_2(t)
\end{bmatrix}.
\end{equation}

As the optimal feedback policy cancels out the non-linearities in the dynamics, \eqref{eq:SMIB:autonomous:sys} is linear and time-invariant. Hence, its stability can be inferred from the eigenvalues of the $A$ matrix, which in this case are $\texttt{eig}(A) = \pm C_2{C_1}^{-1}$.

The positive-real eigenvalue implies the autonomous system will be unstable, see Figure~\ref{fig:SMIB}. This example highlights the challenges of the indirect approach: even for problems that are sufficiently simple to obtain the analytic optimal feedback law, such a policy may not be suitable for closed-loop implementation. 

\subsubsection{Solution via the direct methods}

\begin{figure} 
    \centering
    \includegraphics[width=\columnwidth]{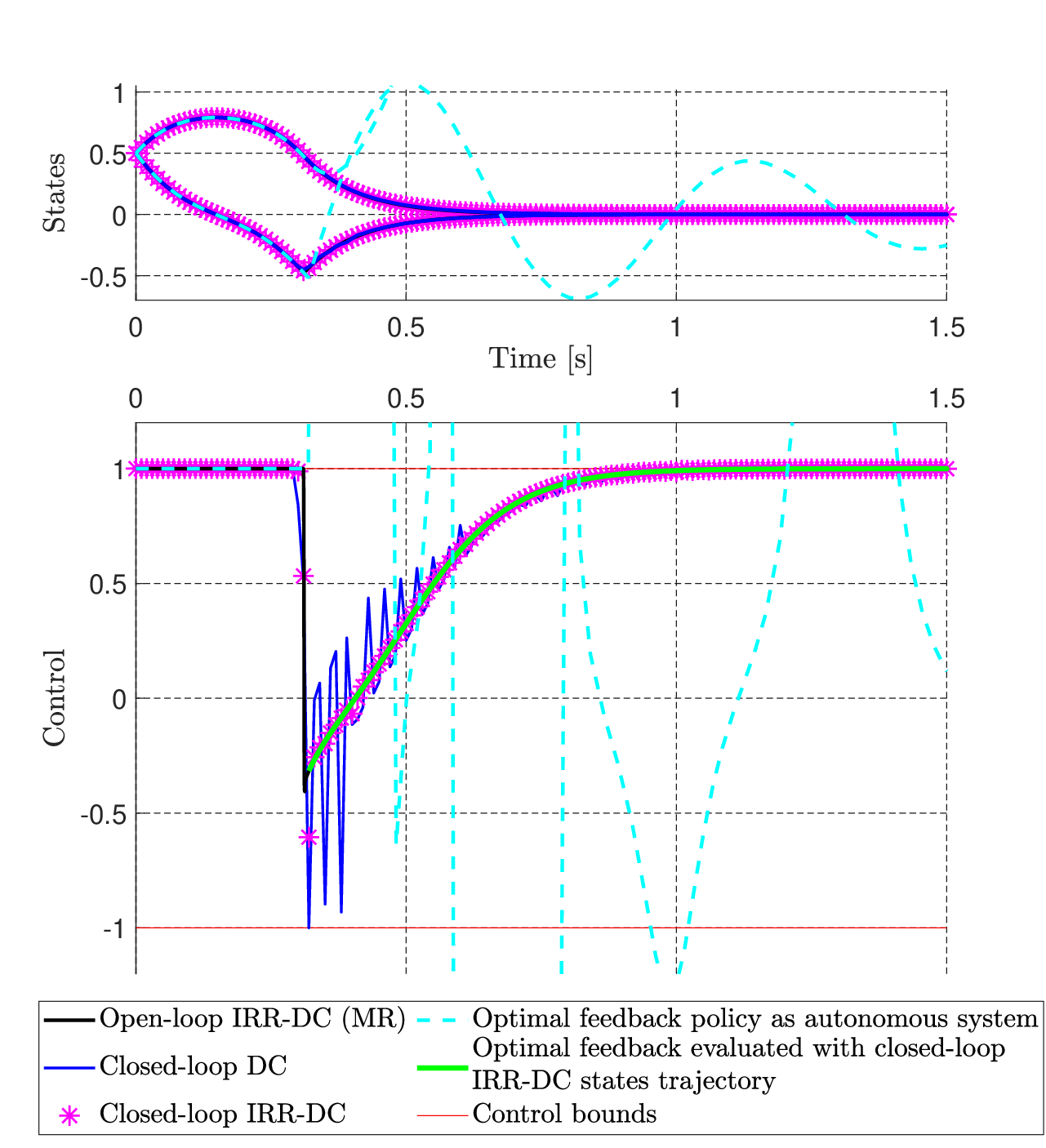}
    \caption{Solutions to the SMIB problem with different methods. The state plot represents $[\frac{x_1}{C_1}, \frac{x_2}{C_2}]^\top$ over time.} 
    \label{fig:SMIB}
\end{figure}

The open-loop and closed-loop solutions using \ac{IRR-DC} are plotted in Figure \ref{fig:SMIB} against the closed-loop DC solution. The solution has a bang-singular structure, \emph{i.e.} after the first switch, $p^*_2~=~0$ throughout. On the singular arc, the optimal feedback policy~\eqref{eq:SMIB:control:feedbackpolicy} evaluated based on the optimal state trajectory obtained by \ac{IRR-DC} is also included for comparison.

\begin{figure}
    \centering
    \includegraphics[width=\columnwidth]{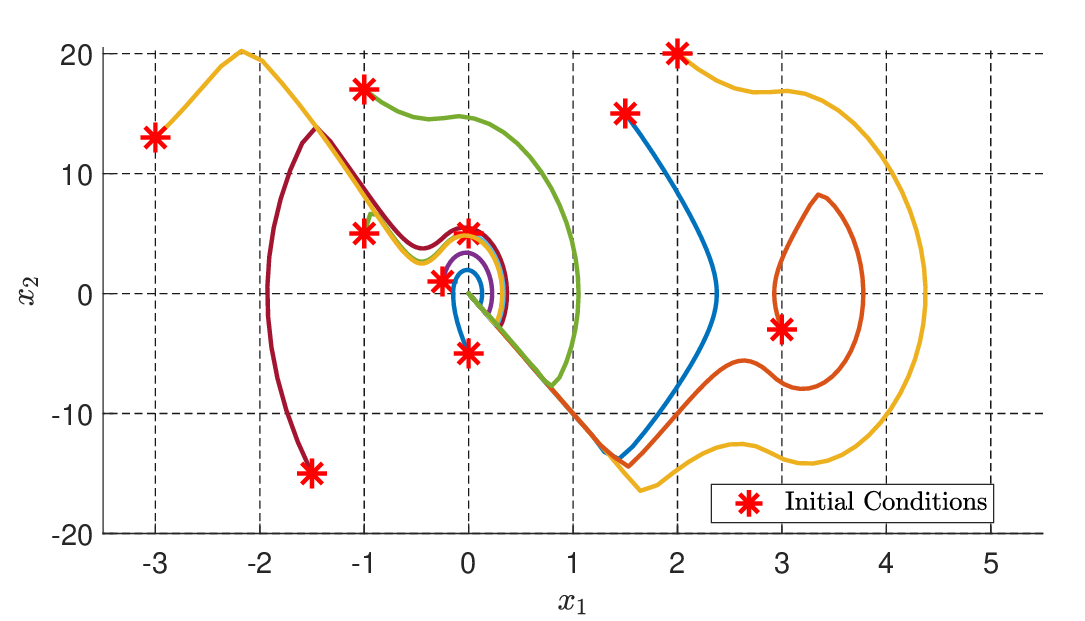}
    \caption{Phase portrait showing the optimal trajectories of the SMIB system for various initial conditions computed using receding horizon implementation of IRR-DC method. }
    \label{fig:phase}
\end{figure}

The closed-loop implementation of IRR-DC provides a control trajectory with the singular arc fluctuation suppressed and agrees with the input computed using the optimal feedback policy. The solutions also pass the Kelley condition \eqref{eq:Singular:optimality}. On closer inspection, the Kelly condition, which guarantees optimality along the singular arc, is only met when $x_1~\geq~-\delta_{ep}$. Since the objective is to drive the system to the origin, \emph{i.e.,}\ $[x_1, x_2]^\top \rightarrow [0,0]^\top$, the optimal trajectory must approach the origin from  $x_1~\geq~-\delta_{ep}$. A similar behaviour occurs when the \ac{OCP} is initialised at different points; see Figure \ref{fig:phase}. Regardless of the initial conditions, the optimal state trajectories obtained via closed-loop \ac{IRR-DC} converge to the singular arc on the fourth quadrant first before arriving at the origin along the singular arc. 

This example demonstrates the benefit of implementing OCP solutions in a closed-loop with EMPC, as opposed to using the analytically derived optimal feedback policy. For solving the OCPs in EMPC, IRR-DC again demonstrates clear advantages over direct collocation with the large fluctuations successfully suppressed.

\section{Conclusions and future work}

This paper examines the closed-loop implementation of optimal control problems involving singular arcs from a practitioner's perspective. We emphasize the importance and the challenges of achieving high-accuracy solutions on the singular arc. Many commonly used methods, such as the indirect approach and the direct collocation method with ad-hoc fixes, require prior knowledge of the existence of singular arcs and/or the structure of the optimal control solution. Such requirements make these methods impractical in closed-loop settings.

Via the examples, we demonstrated the advantages of solving optimal control problems with singular arcs using the direct approach of integrated residual methods. These benefits are particularly evident in closed-loop implementations with economic model predictive control, including the suppression of fluctuations, close agreement with analytical solutions, and the ability to achieve these outcomes without needing prior knowledge of the singular arc to apply special treatments. 

Future work will investigate the interplay between open-loop OCP solutions and closed-loop performance. Specifically, how different choices of model predictive control configurations (\emph{e.g.,}\ mesh size for DC transcription, time step, horizon length) affect the closed-loop cost and the suppression of singular arc fluctuations present in open-loop solutions.

\printbibliography

@Article{Grune2014,
  author          = {Faulwasser, Timm and Gr{\"{u}}ne, Lars and M{\"{u}}ller, Matthias A},
  journal         = {Found. Trends{\textregistered} Syst. Control},
  title           = {{Economic Nonlinear Model Predictive Control}},
  year            = {2018},
  issn            = {2325-6818},
  number          = {1},
  pages           = {224--409},
  volume          = {5},
  doi             = {10.1561/2600000014},
}

@book{betts2010practical,
  title={Practical Methods for Optimal Control and Estimation Using Nonlinear Programming: Second Edition},
  author={Betts, J. T.},
  isbn={9780898718577},
  lccn={2009025106},
  series={Advances in Design and Control},
  year={2010},
  publisher={Society for Industrial and Applied Mathematics}
}

@article{XavierMarceloA.2015Lbcc,
  title={Optimal fast charging of lithium ion batteries: {B}etween model-based and data-driven methods},
  author={Tucker, George and Drummond, Ross and Duncan, Stephen R},
  journal={Journal of The Electrochemical Society},
  volume={170},
  number={12},
  pages={120508},
  year={2023},
  publisher={IOP Publishing}
}

@article{JacobsonD.1970Coos,
issn = {0018-9286},
journal = {IEEE Transactions on Automatic Control},
pages = {67--73},
volume = {15},
publisher = {IEEE},
number = {1},
year = {1970},
title = {Computation of optimal singular controls},
copyright = {Copyright 2004 Elsevier B.V., All rights reserved.},
language = {eng},
author = {Jacobson, D. and Gershwin, S. and Lele, M.},
}

@article{rao2009survey,
  title={A survey of numerical methods for optimal control},
  author={Rao, A. V.},
  journal={Advances in the Astronautical Sciences},
  volume={135},
  number={1},
  pages={497--528},
  year={2009},
  publisher={Univelt, Inc.}
}

@article{wachter2006implementation,
  title={On the implementation of an interior-point filter line-search algorithm for large-scale nonlinear programming},
  author={W{\"a}chter, Andreas and Biegler, Lorenz T},
  journal={Mathematical Programming},
  volume={106},
  number={1},
  pages={25--57},
  year={2006},
  publisher={Springer}
}

@article{Shen1999,
author = {Shen, Haijun and Tsiotras, Panagiotis},
title = {Time-Optimal Control of Axisymmetric Rigid Spacecraft Using Two Controls},
journal = {Journal of Guidance, Control, and Dynamics},
volume = {22},
number = {5},
pages = {682-694},
year = {1999},
doi = {10.2514/2.4436},

URL = { 
    
        https://doi.org/10.2514/2.4436
    
    

},
eprint = { 
    
        https://doi.org/10.2514/2.4436
    
    

}

}

@inproceedings{ICLOCS2,
  title={{ICLOCS2}: {T}ry this optimal control problem solver before you try the rest},
  author={Nie, Yuanbo and Faqir, Omar and Kerrigan, Eric C},
  booktitle={Procs of the 12$^\text{th}$ UKACC International Conference on Control (CONTROL)},
  pages={336--336},
  year={2018},
  organization={IEEE}
}

@article{neuenhofen2018dynamic,
  title={Dynamic Optimization with Convergence Guarantees},
  author={Neuenhofen, Martin P and Kerrigan, Eric C},
  journal={arXiv preprint arXiv:1810.04059},
  year={2018}
}

@book{kirk2004optimal,
  title={Optimal Control Theory: An Introduction},
  author={Kirk, D.E.},
  isbn={9780486434841},
  lccn={74141636},
  series={Dover Books on Electrical Engineering Series},
  url={https://books.google.co.uk/books?id=fCh2SAtWIdwC},
  year={2004},
  publisher={Dover Publications}
}

@book{bryson1975applied,
  title={Applied Optimal Control: Optimization, Estimation and Control},
  author={Bryson, A.E.},
  isbn={9780891162285},
  lccn={75016114},
  series={Halsted Press book'},
  url={https://books.google.co.uk/books?id=P4TKxn7qW5kC},
  year={1975},
  publisher={Taylor \& Francis}
}

@article{KELLEYHENRYJ1964Asvt,
issn = {0001-1452},
journal = {AIAA journal},
pages = {1380--1382},
volume = {2},
number = {8},
year = {1964},
title = {A second variation test for singular extremals},
copyright = {Copyright 2016 Elsevier B.V., All rights reserved.},
language = {eng},
author = {Kelley, HENRY J},
}

@ARTICLE{nie2020efficient,
  author={Nie, Yuanbo and Kerrigan, Eric C.},
  journal={IEEE Control Systems Letters}, 
  title={Efficient and More Accurate Representation of Solution Trajectories in Numerical Optimal Control}, 
  year={2020},
  volume={4},
  number={1},
  pages={61-66},
  doi={10.1109/LCSYS.2019.2921704}}

@article{jones2024modelpredictivebangbangcontroller,
  title={Model Predictive Bang-Bang Controller Synthesis via Approximate Value Functions},
  author={Jones, Morgan and Nie, Yuanbo and Peet, Matthew M},
  journal={IFAC-PapersOnLine},
  volume={58},
  number={17},
  pages={127--132},
  year={2024},
  publisher={Elsevier}
}

@article{ALY01041973,
author = {G. M. Aly and W. C. Chan},
title = {Application of a modified quasilinearization technique to totally singular optimal control problems† },
journal = {International Journal of Control},
volume = {17},
number = {4},
pages = {809--815},
year = {1973},
publisher = {Taylor \& Francis},
doi = {10.1080/00207177308932423},


URL = { 
    
        https://doi.org/10.1080/00207177308932423
    
    

},
eprint = { 
    
        https://doi.org/10.1080/00207177308932423
    
    

}

}

@article{ANDRESMARTINEZ2020106923,
title = {An indirect approach for singular optimal control problems},
journal = {Computers \& Chemical Engineering},
volume = {139},
pages = {106923},
year = {2020},
issn = {0098-1354},
doi = {https://doi.org/10.1016/j.compchemeng.2020.106923},
url = {https://www.sciencedirect.com/science/article/pii/S0098135420303458},
author = {Oswaldo Andrés-Martínez and Lorenz T. Biegler and Antonio Flores-Tlacuahuac}
}

@Article{Pager2022,
author={Pager, Elisha R.
and Rao, Anil V.},
title={Method for solving bang-bang and singular optimal control problems using adaptive {R}adau collocation},
journal={Computational Optimization and Applications},
year={2022},
month={Apr},
day={01},
volume={81},
number={3},
pages={857-887},
issn={1573-2894},
doi={10.1007/s10589-022-00350-6},
url={https://doi.org/10.1007/s10589-022-00350-6}
}

@article{ALY01111978,
author = {G. M. Aly},
title = {The computation of optimal singular control},
journal = {International Journal of Control},
volume = {28},
number = {5},
pages = {681--688},
year = {1978},
publisher = {Taylor \& Francis},
doi = {10.1080/00207177808922489},


URL = { 
    
        https://doi.org/10.1080/00207177808922489
    
    

},
eprint = { 
    
        https://doi.org/10.1080/00207177808922489
    
    

}

}

@phdthesis{pagerThesis,
    author = {Elisha Rose Pager},
    title = {Computational methods for solving nonsmooth and singular optimal control problems with aerospace applications},
    school = {University of Florida},
    year = {2022}
}

@article{Nie_2023,
   title={Solving Dynamic Optimization Problems to a Specified Accuracy: {A}n Alternating Approach Using Integrated Residuals},
   volume={68},
   ISSN={2334-3303},
   url={http://dx.doi.org/10.1109/TAC.2022.3144131},
   DOI={10.1109/tac.2022.3144131},
   number={1},
   journal={IEEE Transactions on Automatic Control},
   publisher={Institute of Electrical and Electronics Engineers (IEEE)},
   author={Nie, Yuanbo and Kerrigan, Eric C.},
   year={2023},
   month=jan, pages={548–555} }

@ARTICLE{ChangChow1998,

  author={Jaewon Chang and Chow, J.H.},

  journal={IEEE Transactions on Power Systems}, 

  title={Time-optimal control of power systems requiring multiple switchings of series capacitors}, 

  year={1998},

  volume={13},

  number={2},

  pages={367-373},

  doi={10.1109/59.667353}}

@ARTICLE{Anderson1972,
  author={Anderson, G.},
  journal={IEEE Transactions on Automatic Control}, 
  title={An indirect numerical method for the solution of a class of optimal control problems with singular arcs}, 
  year={1972},
  volume={17},
  number={3},
  pages={363-365},
  doi={10.1109/TAC.1972.1099989}}

@Article{Powers1980ER,
author={Powers, W. F.},
title={On the order of singular optimal control problems},
journal={Journal of Optimization Theory and Applications},
year={1980},
month={Dec},
day={01},
volume={32},
number={4},
pages={479-489},
issn={1573-2878},
doi={10.1007/BF00934035},
url={https://doi.org/10.1007/BF00934035}
}

@inproceedings{neuenhofen2023numerical,
  title={Numerical Comparison of Collocation vs Quadrature Penalty Methods},
  author={Neuenhofen, Martin P and Kerrigan, Eric C and Nie, Yuanbo},
  booktitle={Procs. of the 62$^\text{nd}$  Conference on Decision and Control (CDC)},
  pages={4285--4290},
  year={2023},
  organization={IEEE}
}

@article{nie2025reliablesolutiondynamicoptimization,
  title={Reliable Solution to Dynamic Optimization Problems using Integrated Residual Regularized Direct Collocation},
  author={Nie, Yuanbo and Kerrigan, Eric C},
  journal={arXiv preprint arXiv:2503.09123},
  year={2025}
}

@misc{harzer2025integrationerrorregularizationdirect,
      title={Integration Error Regularization in Direct Optimal Control using Embedded Runge Kutta Methods}, 
      author={Jakob Harzer and Jochem De Schutter and Moritz Diehl},
      year={2025},
      eprint={2503.12621},
      archivePrefix={arXiv},
      primaryClass={math.OC},
      url={https://arxiv.org/abs/2503.12621}, 
}

@inproceedings{heidrich2021investigation,
  title={Investigation of control regularization functions in bang-bang/singular optimal control problems},
  author={Heidrich, Casey R and Sparapany, Michael J and Grant, Michael J},
  booktitle={AIAA Scitech 2021 Forum},
  pages={2037},
  year={2021}
}

\end{document}